\documentclass[12pt]{article}
\usepackage[dvips]{graphicx}
\DeclareGraphicsExtensions{ps,eps}
\makeatletter
\def\section{\@startsection{section}{1}{\z@}{-1.5ex plus -.5ex
minus -2.ex}{1ex plus .2ex}{\large\bf}}

\input{amssym.def}
\input{amssym}

\def\@thmcounterstep{}
\long\def\@makecaption#1#2{\vskip 10pt \setbox\@tempboxa\hbox{#1.#2}
\ifdim \wd\@tempboxa >\hsize
#1.#2\par 
\else99
\hbox to\hsize{\hfil\box\@tempboxa\hfil}
\fi}
\def\ps@headings{
\def\@oddhead{\footnotesize\rm\hfill\runninghead\hfill}
\def\@evenhead{\@oddhead}
\def\@oddfoot{\rm\hfill\thepage\hfill}\def\@evenfoot{\@oddfoot}}
\newtheorem{Theorem}{Theorem}[section]

\newtheorem{Definition}[Theorem]{Definition}

\makeatother
\title
{
Optimal Transport, Convection, Magnetic Relaxation
and Generalized Boussinesq equations
}
\def\runninghead{\quad 
Convection and optimal transport
}
\author{
{\em Yann Brenier}\thanks{CNRS, Universit\'e de Nice (FR 2800 W. D\"oblin),
Institut Universitaire de France, visiting
Universit\"at Bonn (Institut f\"ur Angewandte Mathematik)
and Universit\"at Wien (Wolfgang Pauli Institut).}}
\date{} 
\begin{document}
\pagestyle{headings}
\flushbottom
\maketitle
\vspace{-10pt}

$$ $$
\section*{Abstract}
$$ $$

We establish a connection between Optimal Transport Theory
(see \cite{Vi} for instance) and 
classical Convection Theory for geophysical flows \cite{Pe}.
Our starting point is the model designed few years ago
by Angenent, Haker and Tannenbaum \cite{AHT} to solve some
Optimal Transport problems.
This model can be seen as a generalization
of the Darcy-Boussinesq equations, which is a degenerate version of
the Navier-Stokes-Boussinesq (NSB) equations.
\\
In a unified framework, we relate different variants of the 
NSB equations (in particular what we call
the generalized Hydrostatic-Boussinesq equations)
to various models involving Optimal Transport
(and the related Monge-Amp\`ere equation \cite{Br,Ca}).
This includes the 2D semi-geostrophic
equations \cite{Ho,CNP,BB,CGP,Lo} and some fully non-linear versions
of the so-called high-field limit 
of the Vlasov-Poisson system \cite{NPS} and of the Keller-Segel
for Chemotaxis
\cite{KS,JL,CMPS}.
\\
Mathematically speaking, we establish some existence theorems
for local smooth, global smooth or global weak
solutions of the different models.
We also justify that the inertia terms
can be rigorously neglected under appropriate
scaling assumptions in the Generalized Navier-Stokes-Boussinesq equations.
\\
Finally, we show how a ``stringy'' generalization of the AHT model
can be related to the magnetic relaxation model studied by
Arnold and Moffatt to obtain stationary solutions of the Euler equations
with prescribed topology (see \cite{AK,Mo,Mo2,Sc,VMI,Ni}).

\section{The Angenent-Haker-Tannenbaum model for Optimal Transport problems}

In this section, we consider the model introduced by 
Angenent, Haker and Tannenbaum \cite{AHT}.
This model was designed in order to 
seek the solutions of some optimal transport problems as equilibrium states
of a suitable dynamical system that could be efficiently solved 
on a computer. The concrete applications 
have been computer vision, image registration and image warping.

\subsection{Optimal transport and rearrangements}

Let us briefly recall some typical results in 
Optimal Transport Theory, such as the polar factorization of maps.
More precisely, let $D$ be the closure of a bounded connected
open set in $R^d$, 
with a boundary of zero $d$-dimensional Lebesgue measure.
Up to a rescaling, we assume the Lebesgue measure of $D$ to be $1$.
Given an  $L^2$ map $y:D\rightarrow R^d$, we call image measure
of the Lebesgue measure on $D$ by $y$ the unique nonnegative
(Borel) measure $\mu$ defined by:
\begin{equation}
\label{image measure}
\int_{R^d}f(x)\mu(dx)=\int_{D}f(y(a))da,
\end{equation}
for all compactly supported continuous functions $f$ on $R^d$.
We have
$$
\int_{R^d}\mu(dx)=1,\;\;\;
\int_{R^d}|x|^2 \mu(dx)=\int_{D}|y(a)|^2 da,
$$
which means that 
$\mu$ belongs to the set $Prob_{2}(R^d)$ of all (Borel)
probability measures $\mu$ on $R^d$ such that
$\int |x|^2 \mu(dx)<\infty$.
In this space, we say that a sequence $\mu_n$
converges tightly to $\mu$ in $Prob_{2}(R^d)$, if:
$$
\int_{R^d} f(x)\mu_n(dx)\rightarrow \int_{R^d} f(x)\mu(dx)
$$
for all continuous function $f$ on $R^d$ such that 
$$
\sup_{x\in R^d}\frac{|f(x)|}{1+|x|^2}\;<+\infty.
$$
Given two $L^2$ maps $y$ and $z$ from $D$ to $R^d$,
we say that they are rearrangement of each other if
they define the same image measure.
When $y$ is a rearrangement of the identity map,
we say, in short, that $y$ is Lebesgue measure preserving.
\\
Next, we define the class of maps with convex potential:
\begin{Definition}   
We say that an $L^2$ map
from $D$ to $R^d$ belongs to the class $C$
of maps with a convex potential if
there is a lower continuous convex function 
$p:R^d\rightarrow ]-\infty,+\infty]$
such that, for Lebesgue almost every point $x\in D$, $p$ is
differentiable at $x$ and its gradient
$\nabla p(x)$ coincides with $y(x)$.
\end{Definition}

Then, we get from \cite{Br}:
\begin{Theorem}   
\label{YB} (Rearrangements with convex potentials)
\\
For each
$L^2$ map $y:D\rightarrow R^d$
there is a unique rearrangement map with a convex potential $y^*\in C$.
The map $y^*$ depends on $y$ only through the associated measure $\mu$
defined by (\ref{image measure}).
\\
In addition, 
the nonlinear operator $y\in L^2\rightarrow y^*\in L^2$ is continuous  
as well as the induced operator
$\mu\in Prob_2(R^d)\rightarrow y^*\in L^2$, 
with respect to the tight convergence.

\end{Theorem}

We get more precise results if $y$ is a non degenerate
map, in the sense that the pre-image
of every Lebesgue negligible set is also negligible:

\begin{Theorem}   
\label{YB2}(Polar factorization of maps \cite{Br})
\\
Let $y$ be a non degenerate $L^2$ map from $D$ to $R^d$.
Then, there is a unique ``polar factorization'' $y=Y\circ X$
where $Y$ belongs to $C$ and $X$ is a Lebesgue measure preserving
map of $D$.
In this decomposition, $Y$ is the unique rearrangement $y^*$ of $y$ in $C$
and $X$ is the unique measure preserving map of $D$ that
minimizes
$\int_D|X(a)-y(a)|^2\;da.$
In addition, $X$ can be written:
\begin{equation}
\label{polar1}
X(a)=(\nabla\Phi)(y(a)),\;\;\;a.e.\;a\in D,
\end{equation}
where $\Phi$ is a convex Lipschitz function defined on $R^d$.
\end{Theorem}

For (much) more results on optimal transport, we refer to Villani's
textbook \cite{Vi}.
The expression ``optimal transport'' comes from the fact that $y^*$,
among all possible rearrangements $y$ of $y^0$, is the unique minimizer of
the ``transportation cost''. 
$$
\int_D |y(x)-x|^2 \;dx,
$$
where $|\cdot|$ denotes the Euclidean norm.
The name ``map with convex potential''
is due to Caffarelli \cite{Ca}. The concept of polar factorization has
been extended to Riemannian manifolds by McCann \cite{Mc}.
Examples of concrete applications of optimal transport techniques to 
natural and computer sciences can be found in \cite{FMMS,HZTA}.

\subsection{The AHT model}

The AHT model is an attempt to get the unique rearrangement $y^*$ 
of $y^0$, with convex potential,
as the equilibrium state at $t=+\infty$ of the following set
of evolution equations:
\begin{equation}
\label{AHT1}
\partial_t  y+(v\cdot\nabla) y=0,
\end{equation}
\begin{equation}
\label{AHT2}
Kv+\nabla p= y,\;\;\;\nabla.v=0,
\end{equation}
where $y=y(t,x)\in R^d$, 
$v=v(t,x)\in R^d$, $p=p(t,x)\in R$
depend on $t\ge 0$ and
$x\in D$, and 
$K$ is a ``dissipative'' operator to be chosen, for instance $K=I$
or $K=-\Delta$.
In these ``AHT'' equations, we denote
the inner product in $R^d$ by $\cdot$ 
and we use notations:
$$
\nabla_i=\frac{\partial}{\partial x_i},\;\;\;
v\cdot\nabla=\sum_{j=1,d}\;v_j\frac{\partial}{\partial x_j},
\;\;\;\Delta=\sum_{j=1,d}\;\frac{\partial^2}{\partial x_j^2}.
$$
The boundary conditions for the AHT system (\ref{AHT1},\ref{AHT2}) are:
\\
i) the initial value of $y$ at $t=0$, $y(0,x)=y^0(x)$,
\\
ii) $v$ is parallel to the boundary $\partial D$ if $K=I$
and $v=0$ along the boundary if $K=-\Delta$.
\\
Notice that neither $p$ nor $v$ need initial conditions.
As a matter of fact, as $K=I$, the second AHT equation (\ref{AHT2})
just corresponds to the ``Helmholz decomposition'' of $y$ as a sum
of a gradient field and a divergence-free field parallel to the
boundary $\partial D$. The field $p$ can be recovered by solving
the Poisson problem:
$$
\Delta p=\nabla.y,
$$
inside $D$ with inhomogeneous Neumann condition $\nabla p\cdot n=y\cdot n$
along the boundary, where $n$ denotes the outward normal. 
Then, we get: $v=y-\nabla p$.
So, we can write:
$v=Py$, where $P$ is a linear singular integral operator bounded in
all $L^p$ space for $1<p<+\infty$, provided that the domain $D$ is
smooth enough.
In the case $K=-\Delta$, in a similar way we can write $v=P_\Delta y$,
where $P_\Delta$ is a linear singular integral operator bounded 
from $L^p$ to the Sobolev space $W^{2,p}$, 
for all $1<p<+\infty$, $D$ being assumed to be smooth.
Thus we can write the AHT system (\ref{AHT1},\ref{AHT2}) 
in a more abstract form:
\begin{equation}
\label{compact AHT}
\partial_t  y+(P_K y\cdot\nabla) y=0,
\end{equation}
with $P_K=P$ if $K=I$ and $P_K=P_\Delta$ if $K=-\Delta$.

\subsection{Expected long time behaviour of the AHT model}
\label{long time}

Let us now explain why the AHT model is expected to
solve the Optimal Transport (or rearrangement) problem, at least for
a large class of data.
First, we observe that equation (\ref{AHT1}) expresses,
at least formally, that,
at each time $t$, $y(t,\cdot)$ is a rearrangement of $y^0$.
Indeed, for any smooth compactly supported function $f$, we get:
$$
\frac{d}{dt}\int_D f(y(t,x)) \;dx=
-\int_D (\nabla f)(y)\cdot(v\cdot\nabla)y \;dx
$$
$$
=-\int_D v\cdot\nabla[f(y)]\;dx,
$$
(using the chain rule)
which is zero, since 
$v$ is divergence free and parallel to $\partial D$ and is
therefore $L^2$ orthogonal to any gradient field.
(Notice that this calculation can be made rigorous provided that $v$ has
enough regularity. According to Ambrosio's recent improvement of
the DiPerna-Lions theory on ODEs \cite{Am,DL}, it is enough that
$v$ belongs to $L^1_{loc}(R_+,BV(D,R^d))$.)
\\
Next, we get the following balance law for the AHT 
\begin{equation}
\label{decay}
\frac{d}{dt}\int_D \frac{1}{2}|y(t,x)-x|^2 \;dx=-\int_D (v\cdot Kv)(t,x)\;dx,
\end{equation}
which implies, at least formally, that
$y$ and $v$
respectively belong to the functional spaces $L^\infty(R_+,L^2(D,R^d))$
and $L^2(R_+,H_K(D))$. Here $H_K(D)$ denotes the Hilbert space of all 
divergence free fields $w(x)\in R^d$ 
for which $\int_D w\cdot Kw \;dx$ is finite
with suitable boundary conditions ($w$ parallel to $\partial D$ as
$K=I$ or $w=0$ on $\partial D$ as $K=-\Delta$).
The formal proof of (\ref{decay}) is as follows:
$$
\frac{d}{dt}\int_D \frac{1}{2}|y-x|^2 \;dx=
-\int_D (y-x)\cdot ((v\cdot\nabla)y)\;dx
$$
(using the first AHT equation (\ref{AHT1}))
$$
=-\int_D v\cdot\nabla(\frac{1}{2}|y-x|^2) \;dx-\int_D v\cdot(y-x)\;dx
$$
$$
=-\int_D v\cdot(y-x)\;dx
$$
(since $v$ is divergence free and parallel to $\partial D$ and is
therefore $L^2$ orthogonal to any gradient field)
$$
=-\int_D v\cdot(Kv+\nabla (p-\frac{1}{2}|x|^2)) \;dx
$$
(using the second AHT equation (\ref{AHT2}))
$$
=-\int_D v\cdot Kv \;dx.
$$
At this point, we can describe the expected long time behaviour of the AHT 
system through the following heuristics.
Since $Kv\cdot v$ is space-time integrable, we first argue
that, as $t\rightarrow +\infty$,
$v$ presumably tends to zero.
Then we expect $y$ to have a definite strong limit $y^\infty$ in $L^2$,
which is, then,  necessarily a rearrangement of $y^0$.
Passing to the limit in the second AHT equation (\ref{AHT2}), we conclude that
$y^\infty$ must be a gradient. Therefore, at the end of the process,
$y^0$ has been rearranged as a map  $y^\infty$ with a potential.
Observe that this potential needs not being convex.
This is obvious in the special case 
when $y^0$ is itself a map with a potential which is not convex.
Indeed, then
$$
y(t,x)=y^0(x),\;\;\;v(t,x)=0,
$$
is a trivial stationary solution 
to the AHT equations (\ref{AHT1},\ref{AHT2}) and we get
$y^\infty=y^0$ as a map with a non convex potential.
So, we need further assumptions on $y^0$ to be convinced that $y^\infty$
has a chance to have a convex potential.
A natural assumption is that $y^0$ is smooth with a
positive jacobian determinant valued in some interval $[r,1/r]$ with
$0<r<+\infty$.
Indeed, for such an initial condition,
the AHT equations have a global smooth solution $y$ (at least in the case
when $K=-\Delta$, as discussed in the next subsection),
with a jacobian determinant that must stay in the same
interval $[r,1/r]$, since $y$ is always a rearrangement of $y^0$.
So if the convergence to $y^\infty$, 
as $t\rightarrow +\infty$, is strong enough, we expect $y^\infty$ to have a
convex potential and, therefore, coincide with
the unique rearrangement of $y^0$ with convex
potential $y^*$ provided by Theorem \ref{YB}.
The results obtained in \cite{AHT} are only partial and leave as an open
question this issue.

\subsection{Wellposedness of the AHT equations}

From the PDE viewpoint, it is crucial to check that the
AHT system (\ref{AHT1},\ref{AHT2}) is wellposed, 
which is done by Angenent, Haker and Tannenbaum in
\cite{AHT}, for a class of dissipative operator $K$ including $K=I$.
Let us briefly discuss the wellposedness issue in the cases $K=I$
and $K=-\Delta$.
\\
For $K=I$, the AHT system is
similar to the inviscid Burgers equation, 
\begin{equation}
\label{burgers}
\partial_t  y+(y\cdot\nabla) y=0,
\end{equation}
since $P_K$ behaves like a
pseudo-differential operator of order zero.
Thus, the local in time existence of smooth solutions for smooth
initial conditions can be obtained from rather standard energy estimates.
It is a challenging and interesting open question 
whether the Lipschitz norm, in space, of such solutions may blow
up in finite time 
(as it would be the case of the inviscid Burgers equation).
In sharp contrast, in the case $K=-\Delta$, 
smooth solutions are clearly global in time.
Indeed, from (\ref{AHT1}), we immediately get
that $|y(t,x)|$ is uniformly bounded by the sup norm of 
$y^0$ that we denote by $M^0$ and suppose, here, to be finite.
Thus, because of (\ref{AHT2}), according to standard elliptic regularity
theory, the $L^\infty_t(W_x^{2,p})$ norm of $v(t,x)$ is
controled by $M^0$ for all finite $p$. Thus, the same is true for the
sup norm of 
$\nabla v(t,x)$. Differentiating (\ref{AHT1}) in $x$, we deduce
that the sup norm of $\nabla y(t,x)$ in $x$
cannot grow, in sup norm, faster than exponentially in $t$ as soon
as $\nabla y^0$ has a finite sup norm.
So, there is no possible blow up of the Lipschitz norm of
both $v$ and $y$ and, therefore, by a standard argument,
smooth solutions must be global in time.
(Notice that the dissipative operator $K=(-\Delta)^{1/2}$, with
appropriate boundary condition, would be borderline to get
such a Lipschitz estimate. In that case $v$ would be a priori only
Log-Lipschitz, just like the Yudovich solutions of the 2D Euler
equations \cite{MP}.)
For more details, we refer to the paper by Angenent, Haker
and Tannenbaum \cite{AHT}, where different kinds of operator $K$
are considered.

\subsection{Interpretation of the AHT system in terms of Convection Theory}

From a Fluid Mechanics viewpoint, the AHT equations look very similar
to the Boussinesq equations for convective flows, in particular
to their Darcy-Boussinesq version. 
A classical model for Convection Theory is provided by
the Navier-Stokes Boussinesq (NSB) equations that we are now going
to review with more details.
Using the Boussinesq approximation, 
the Navier-Stokes equations for an inhomogeneous
incompressible fluid subject to gravity along the $x_d$ direction read:
\begin{equation}
\label{NSB1}
\rho_0(\partial_{t}v+(v\cdot\nabla)v)-\nu\Delta v+\nabla p=y,
\;\;\;\nabla\cdot v=0,
\end{equation}
\begin{equation}
\label{NSB2}
\partial_{t}y+(v\cdot\nabla)y=0.
\end{equation}
Here, $v=v(t,x)\in R^d$ is the velocity field, $p=p(t,x)\in R$ the
pressure field, $\rho_0>0$ the average density of the fluid,
$\nu$ the (constant) viscosity of the fluid, while
$y$ has only one component in the ``vertical'' direction $x_d$,
which is $-g\;\theta(t,x)$, where $g$ is the gravity constant
and $\theta(t,x)$ is the difference between the
density of the fluid at $(t,x)$ and the averaged density $\rho_0$ of the
fluid.
(Usually in Convection Theory, 
a diffusion term is added to the advection equation for $y$
\cite{ID}.)
We recall that the Boussinesq approximation amounts to consider a variable
density incompressible fluid for which the density variations are sufficiently
small to be neglected in all terms except the gravity force.
This approximation is widely used for
ocean and atmosphere modelling \cite{Pe}. (To the best of our knowledge
the justification of this approximation is still an open
problem in mathematical Fluid Dynamics,
mostly because of our rather poor knowledge of the Navier-Stokes
equations for inhomogeneous flows, see discussions
in \cite{Li,Ma} for instance.)
By neglecting the inertia term, or equivalently by setting $\rho_0=0$
in the NSB (Navier-Stokes-Boussinesq) equations, 
we get the simpler Stokes-Boussinesq (SB) 
(related to large-Prandtl-number Convection Theory as in \cite{DOR,Wa}, 
for instance). If, in addition,
the diffusion term $-\nu\Delta v$ is replaced by a friction
term such as $v$, we get the Darcy-Boussinesq (DB) model.
We immediately see that both the SB and the DB equations are just
particular cases of the
AHT model (\ref{AHT1},\ref{AHT2}), for which the vector valued function $y$
has only one component along the $x_d$ axis. Indeed, the DB and SB models
then respectively correspond to the choice $K=I$ and $K=-\Delta$ in the
second AHT equation (\ref{AHT2}).
According to the discussion made in subsection \ref{long time},
we expect, for the AHT model,
the $y(t,x)$ to converge, as $t\rightarrow +\infty$,
to a map $y^\infty(x)$  with, hopefully, a convex potential.
In the particular case of the convective DB and SB models, 
$y(t,x)$ has only one component in the $x_d$ direction, namely
$-g\;\theta(t,x)$. Interestingly enough,
the convergence toward a map with convex potential,
$exactly$ means, for the DB and SB models, that
the density field tends to a density ``profile''
$\rho^\infty(x_d)$, depending
only on the vertical coordinate $x_d$, and monotonically decreasing.
This clearly corresponds, in terms of Convection Theory, to
a ``stable hydrostatic equilibrium''. 
Notice that a similar discussion can be found in Moffatt's paper
\cite{Mo} (section 2) as a prelude to his Magnetic Relaxation model
that we will consider at the end of the present paper.

\section{Generalized Navier-Stokes-Boussinesq equations}

The interpretation of the AHT model in terms of Convection Theory 
suggests the following ``GNSB'' generalization of the 
NSB (Navier-Stokes-Boussinesq) equations:
\begin{equation}
\label{GNSB1}
\epsilon(\partial_{t}v+(v\cdot\nabla)v)+K v+\nabla p=F(x,y),
\;\;\;\nabla\cdot v=0,
\end{equation}
\begin{equation}
\label{GNSB2}
\partial_{t}y+(v\cdot\nabla)y=G(x,y),
\end{equation}
where $y=y(t,x)\in R^m$ is a vector-valued function ($m\ge 1$,
in practice $m=d$ or $m=2d$ for the models discussed below),
$F$ and $G$ are given smooth functions with bounded derivatives up to
second order, respectively defined on $R^m$ and $R^m\times R^d$,
$\epsilon>0$ is a scaling factor introduced to single out the inertia term,
and $K$ is a linear dissipative operator.
Depending on the applications in view,
only the following cases will be considered:
$K=0$ (no dissipation), $Kv=v$ (linear friction), $Kv=-\Delta v$ (viscosity).

\subsection{Existence theory for the GNSB equations}

For simplicity, we consider in this subsection the domain $D$ to be the unit
periodic cube $T^d=R^d/Z^d$, in order to avoid technicalities due to
spatial boundary conditions.
For the three possible choices of the dissipative operator $K$
\begin{equation}
\label{K}
Kv=0,\;\;\;
Kv=v,\;\;\;
Kv=-\Delta v,
\end{equation}
the existence and uniqueness of a local in time
smooth solution $(y,v)$ of the GNSB equations
(\ref{GNSB1},\ref{GNSB2}), for each smooth initial initial condition 
$(y^0,v^0)$ given on the torus $T^d$, 
follow from standard theory on Euler and Navier-Stokes
equations (for which we refer to \cite{Li}).
\\
We say that $(y,v)$ is a weak solution if:
\\
1) $y(t,x)$ and $v(t,x)$ depends continuously on $t$ with values in 
$L^2(D,R^d)$ 
(with respect to the weak topology of $L^2$);
\\
2) For all smooth time dependent vector fields
$w(t,x)$, $z(t,x)$, with $\nabla\cdot w=0$, we have:
\begin{equation}
\label{GNSB1-weak}
\frac{d}{dt}\int v\cdot w\;dx=
\int [\epsilon v\cdot(\partial_t+(v\cdot\nabla))w-Kv\cdot w+F(x,y)\cdot w]dx,
\end{equation}
\begin{equation}
\label{GNSB2-weak}
\frac{d}{dt}\int y\cdot z\;dx=
\int y\cdot(\partial_t+(v\cdot\nabla))z+G(x,y)\cdot z)dx,
\end{equation}
3) The following energy inequality holds true:
\begin{equation}
\label{energy-weak}
\frac{1}{2}\frac{d}{dt}
\int(\;\epsilon |v|^2+|y|^2\;)dx+\int\;Kv\cdot v\;dx
\le \int [F(x,y)\cdot v+G(x,y)\cdot y]\;dx.
\end{equation}
When $K=-\Delta$, 
the existence of global weak solutions for the GNSB
equations follows from standard arguments \`a la Leray combined with
the DiPerna-Lions theory on ODEs \cite{Li,DL}. They are unique
in 2 space dimensions. In sharp contrast, as $K=I$ or $K=0$, nothing 
can be said about global weak solutions.
\\
Concerning global smooth solutions, the existence theory is quite
challenging, even in 2 space dimensions.
Recently, Chae, Hou and Li \cite{Ch,HL} have proven that
Navier-Stokes Boussinesq equations (\ref{NSB1},\ref{NSB2})
(just called Boussinesq equations in these papers)
have global smooth solutions when $d=2$ and $K=-\Delta$.
The same result can be readily extended to the GNSB equations 
(\ref{GNSB1},\ref{GNSB2}),
essentially because we assume the right-hand sides
$F(x,y)$ and $G(x,y)$
of each equation to be smooth functions of $y$ and $v$
with bounded derivatives up to order two.
(Indeed, these assumptions are enough for a straghtforward adaptation 
of the proof of Theorem 1.1
in Chae's paper, through estimates $(2.1\cdot\cdot\cdot 17)$
in \cite{Ch}. Some constants involved in these estimates have just to
be modified to take into account 
the Lipschitz constants of $F$ and $G$.)
\\
So, we can summarize all these results in the following Theorem,
which is nothing but a straightforward adaptation of known results:

\begin{Theorem}   
\label{chae GNSB}
Assume the dissipative operator $K$ to be of type (\ref{K}).
Then, the generalized Navier-Stokes Boussinesq equations 
(\ref{GNSB1},\ref{GNSB2}) admit, for any smooth initial condition,
a unique local smooth solution.
\\
If $K=-\Delta$, the GNSB equations admit at least
a global weak solution $(y,v)$ 
(in the sense 
(\ref{GNSB1-weak},\ref{GNSB2-weak},\ref{energy-weak}))
for any initial condition $(y^0,v^0)$ in $L^2$.
If $d=2$, these weak solutions are unique. Furthermore, still for $d=2$,
the solutions are globally smooth for smooth initial conditions.
\end{Theorem}

\subsection{Zero inertia limit of the GNSB equations}

By zero inertia limit of the GNSB, we mean the formal limit obtained
by dropping the scaling factor $\epsilon$ in front of the inertia terms
in (\ref{GNSB1},\ref{GNSB2}).
Namely, in Eulerian coordinates,
\begin{equation}
\label{HF}
\partial_{t}y+(v\cdot\nabla)y=G(x,y),
\end{equation}
$$
K v+\nabla p=F(x,y),
\;\;\;\nabla\cdot v=0.
$$
We are able to make a rigorous derivation of
the zero inertia limit when $K$ is strictly
dissipative ($K=0$ being so excluded):

\begin{Theorem}   
\label{high field limit}
Assume the dissipation operator $K$ to be coercive in $L^2$,
namely $K\ge \alpha$
for some constant $\alpha>0$.
Then the zero inertia equations (\ref{HF})
admit, for any smooth initial condition,
a local smooth solution, which is global if $d=2$ and $K=-\Delta$.
This solution can be obtained as the limit,
as $\epsilon$ goes to zero, of the weak
solutions of the GNSB equations (\ref{GNSB1},\ref{GNSB2}),
with the same initial condition.
\end{Theorem}

For the convergence, we use a simple energy method. 
Namely, given a weak solution $(y',v')$ to the GNSB equations 
(\ref{GNSB1-weak},\ref{GNSB2-weak},\ref{energy-weak}) 
and a solution $(y,v)$ of the HF equations, with same
initial conditions $(y^0,v^0)$,
we introduce 
\begin{equation}
\label{relative energy}
e(t)=
\int_{T^d}(\;\epsilon |v'|^2+|y-y'|^2\;)\;dx
\end{equation}
and try to get an estimate of form:
\begin{equation}
\label{gronwall}
\frac{d}{dt}(e(t)+O(\epsilon))+\frac{1}{2}\int_{T^d}K(v-v')\cdot(v-v')\;dx
\le (e(t)+O(\epsilon))c,
\end{equation}
where $c$ depends only on the limit $(y,v)$,
for any fixed finite time interval $[0,T]$ on which $(y,v)$ is smooth.
From this estimate (\ref{gronwall}), we immediately get that
$y-y'$ and $v-v'$ are of order $O(\sqrt\epsilon)$  in, respectively,
$L^\infty([0,T],L^2(T^d))$ and $L^2([0,T],L^2(T^d))$,
using the coercivity of $K$ ($K\ge \alpha$ for some $\alpha>0$).
So, we are left with proving (\ref{gronwall}).
Notice first that, from equations (\ref{energy-weak}) and (\ref{HF}),
the following energy balances hold true:
$$
\frac{1}{2}\frac{d}{dt}
\int(\;\epsilon |v'|^2+|y'|^2\;)dx+\int\;Kv'\cdot v'\;dx
\le\int [F(x,y')\cdot v'+G(x,y')\cdot y']\;dx
$$
$$
\frac{1}{2}\frac{d}{dt}
\int\;|y|^2\;dx+\int\;Kv\cdot v\;dx
=\int [F(x,y)\cdot v+G(x,y)\cdot y]\;dx.
$$
Since $(y',v')$ and $(y,v)$ are respectively supposed to be 
a weak solution of 
the GNSB equations and a smooth solution of the zero inertia limit (\ref{HF}),
we also get
$$
-\frac{d}{dt}
\int\;y\cdot y'\;dx
$$
$$
=\int\;[((v\cdot\nabla)y)\cdot y'-((v'\cdot\nabla)y)\cdot y']dx
-\int [G(x,y)\cdot y'+G(x,y')\cdot y]dx
$$
$$
=\int\;((v\cdot\nabla)y-(v'\cdot\nabla)y)\cdot (y'-y)\;dx
$$
$$
+\int [(G(x,y)-G(x,y'))\cdot(y-y')-G(x,y')\cdot y'-G(x,y)\cdot y]dx
$$
(using that $\int (w\cdot\nabla)y)\cdot y)\;dx=0$ for both $w=v$
and $w=v'$)
$$
\le c_1\int\;( |v-v'||y-y'|+|y-y'|^2)\;dx
-\int [G(x,y')\cdot y'+G(x,y)\cdot y]dx
$$
where $c_1$ depends on the Lipschitz constants of $G$ (as a function of
$y$ and $v$) and $v$ (as a function of $x$).
Thus, adding up these three equalities, we get
by definition (\ref{relative energy}):
$$
\frac{d}{dt}e(t)+\int\;(Kv'\cdot v'+Kv\cdot v)\;dx
$$
$$
\le c_1\int\;( |v-v'||y-y'|+|y-y'|^2)\;dx
+\int [F(x,y')\cdot v'+F(x,y)\cdot v]\;dx.
$$
From (\ref{GNSB1-weak}) and (\ref{HF}), we also get:
$$
\int F(x,y)\cdot v'\;dx=\int (Kv+\nabla p)\cdot v'\;dx=\int Kv\cdot v'\;dx
$$
and
$$
\int F(x,y')\cdot v\;dx=\epsilon \frac{d}{dt}\int v'\cdot v\;dx
-\epsilon\int v'\cdot(\partial_t+v'\cdot\nabla)v\;dx
+\int Kv'\cdot v\;dx
$$
$$
=\sqrt{\epsilon} (\frac{d}{dt}r_1(t)+r_2(t))+e_2(t)+\int Kv'\cdot v\;dx,
$$
where
$$
r_1(t)=\sqrt{\epsilon} \int v'\cdot v\;dx,
$$
$$
r_2(t)=-\sqrt{\epsilon} \int v'\cdot \partial_t v\;dx
$$
and
$$
e_2(t)=-\epsilon \int v'\cdot(v'\cdot\nabla)v\;dx.
$$
Notice that $r_1^2$, $r_2^2$ and $|e_2|$ are bounded by $c_2 e(t)$
(by definition (\ref{relative energy})),
where $c_2$ depends on the Lipschitz constant of $v$ as a function of
both $t$ and $x$.
So, we have obtained:
$$
\frac{d}{dt}(e(t)-\sqrt{\epsilon}r_1(t))
+\int\;K(v'-v)\cdot(v'-v)\;dx
$$
$$
\le c_1\int\;( |v-v'||y-y'|+|y-y'|^2)\;dx+\sqrt{\epsilon}r_2(t)+c_2 e(t).
$$
Using the coercivity of $K$ ($K\ge \alpha>0$)
and definition (\ref{relative energy}), we find $c_3$
depending on the Lipschitz constants of $F$, $G$ and $v$ such that:
$$
c_1\int\;( |v-v'||y-y'|+|y-y'|^2)\;dx \le 
c_3 e(t)+\frac{1}{2}\int\;K(v'-v)\cdot(v'-v)\;dx.
$$
This leads to:
$$
\frac{d}{dt}(e(t)-\sqrt{\epsilon}r_1(t))
+\frac{1}{2}\int\;K(v'-v)\cdot(v'-v)\;dx
\le (c_2+c_3)e(t)+\sqrt{\epsilon}r_2(t),
$$
which leads to a differential inequality of the desired type,
namely (\ref{gronwall}),
since $r_1^2+r_2^2\le 2 c_2 e$.
Thus, the proof of Theorem \ref{high field limit} is now achieved.

\subsection{A Mechanical interpretation of the GNSB equations}
\label{mechanical}

In this subsection, we provide a mechanical interpretation 
of the GNSB equations and their zero inertia limit.
For this purpose,
it is worth considering the
GNSB equations (\ref{GNSB1},\ref{GNSB2}) in ``Lagrangian coordinates''.
Assuming the vector field $v$ to be smooth enough,
denoting by $a\in D$ the
position of a fluid parcel at time $t=0$, we can recover
its position $X(t,a)\in D$
at later time $t\ge 0$  by solving the ODE
\begin{equation}
\label{ode}
\partial_t X(t,a)=v(t,X(t,a)),\;\;\;X(0,a)=a,\;\;\;\forall a\in D.
\end{equation}
Notice that, for each $t$, $a\in D\rightarrow X(t,a)\in D$
is a measure preserving map
as a consequence of the fact that $v$ is a smooth divergence-free
vector field parallel to $\partial D$.
Let us also introduce:
\begin{equation}
\label{Y}
Y(t,a)=y(t,X(t,a))\in R^m.
\end{equation}
Then, 
the GNSB equations (\ref{GNSB1},\ref{GNSB2}) read in Lagrangian coordinates:
\begin{equation}
\label{GNSB-lagrangian}
\epsilon\; \partial_{tt}X(t,a)+(Kv)(t,X(t,a))+
(\nabla p)(t,X(t,a))=F(X(t,a),Y(t,a)),
\end{equation}
$$
\partial_{t}Y(t,a)=G(X(t,a),Y(t,a)),
$$
where $a\in D\rightarrow X(t,a)\in D$ is Lebesgue measure preserving.
Let us now provide a possible mechanical interpretation.
We model the atmosphere
(or the ocean)
as a continuous distribution of infinitesimal rigid balloons
floating inside $D$, each of them having position $X(t,a)$ at time
$t$, with $X(0,a)=a$, and being 
attached with probability $\lambda(a)\ge 0$ to an
an anchor with an elastic cable.
Of course, to be a realistic model with real balloons, $\lambda$
should be a discrete probability distribution concentrated on
a finite collection of points
and the corresponding balloons should
have a finite extension!
Let us rather assume, for mathematical simplicity, that the balloons
are just points and that $\lambda$ is a smooth nonnegative 
density function on $D$ with unit mass.
The cable corresponding to the balloon labelled by $a$ 
is modelled by a (possibly non Hookean) spring with 
restoring force $-k(\xi,a)=k(-\xi,a)\in R^d$ where $\xi\in R^d$ is
the elongation of the spring. Notice that $k$ may depend on $a$.
The location of the anchor attached to the balloon labelled by $a$
is not necessarily fixed and denoted
by $Y(t,a)\in R^d$. (We may also think of an aircraft, or a boat,
or any kind of
carrier instead of an anchor.)
Notice that we do not require the anchor to
be located in $D$.
Neglecting any interaction between the fluid and both the anchors
and the springs (which may not be very realistic),
we obtain the following dynamical equation for each balloon
\begin{equation}
\label{buoys epsilon}
\epsilon\partial_{tt}X(t,a)+(Kv)(t,X(t,a))+(\nabla p)(t,X(t,a))=
\end{equation}
$$
-\lambda(a) k(X(t,a)-Y(t,a),a).
$$
(Observe that as $\lambda(a)=0$,
the corresponding carrier $Y(t,a)$ is just fictitious!)
Let us consider the special case when the speed of each anchor is constant 
and given by:
\begin{equation}
\label{carrier 1}
\partial_t Y(t,a)=W(a)\in R^2.
\end{equation}
Implicitly define a field $y=(\tilde{y},\hat{y})=y(t,x)\in R^d\times D$
by setting 
$$
\tilde y(t,X(t,a))=Y(t,a),\;\;\;\hat y(t,X(t,a))=a
$$
(remember that $a\rightarrow X(t,a)$ is supposed to be a diffeomorphism).
Noticing that
$$
((\partial_t+v\cdot\nabla)y)(t,Y(t,a))=\partial_t[y(t,X(t,a))]
=(\partial_t Y(t,a),0)
$$
$$
=(W(a),0)
=(W(\hat y(t,X(t,a))),0)
$$
and going back to Eulerian coordinates,
we recover the GNSB equations (\ref{GNSB1},\ref{GNSB2}) in the particular case:
\begin{equation}
\label{model 1}
F(x,y)=-\lambda(\hat y)k(x-\tilde y,\tilde y),\;\;\;
G(x,y)=(W(\hat y),0).
\end{equation}
(Notice that assuming $\lambda$ and the restoring force $k$ 
to have bounded derivatives up
to order 2 is not very realistic!
These assumptions are clearly made for mathematical convenience.)
We may also consider the following variant of this mechanical model.
Instead of prescribing their velocity by 
(\ref{carrier 1}), we may assume that the carriers
are driven by a friction-dominated
retroaction of type:
$$
\eta\partial_{tt}Y(t,a)+\partial_t Y(t,a)=-\mu(a)k(Y(t,a)-X(t,a),a),
$$
where $\mu$ is a given nonnegative function.
Dropping the inertia term ($\eta=0$) leads to the following law:
\begin{equation}
\label{carrier 2}
\partial_t Y(t,a)=-\mu(a)k(Y(t,a)-X(t,a),a),
\end{equation}
In that case, we get (keeping unchanged the definitions of $Y$ and $y$)
again the GNSB equations (\ref{GNSB1},\ref{GNSB2}) with $F$ and $G$ given by:
\begin{equation}
\label{model 2}
F(x,y)=-\lambda(\hat y)k(x-\tilde y,\hat y),\;\;\;
G(x,y)=(-\mu(\hat y)k(\tilde y-x,\hat y),0).
\end{equation}
Let us finally consider a second variant where the Coriolis force
is added to the model (rotating ocean or atmosphere). 
Neglecting the vertical extension, so that $d=2$,
and assuming the rotation vector to be perpendicular to the ocean
and of unit length,
the Coriolis force $Jv=(-v_2,v_1)$
is completely absorbed by the pressure term and the fluid parcels
are not sensitive to it. (Indeed, $v=(v_1,v_2)$
being divergence free, $Jv$ is a gradient and can be removed from
the dynamical equation.) However we may think that the carriers
are still sensitive to the Coriolis force.
Thus, we get for them
$$
\eta\partial_{tt}Y(t,a)+J\partial_t Y(t,a)=-\mu(a)k(Y(t,a)-X(t,a),a),
$$
instead of (\ref{carrier 2}), neglecting a possible friction term.
Neglecting the inertia term ($\eta=0$) leads to the balance equation:
\begin{equation}
\label{carrier 3}
\partial_t Y(t,a)=J\mu(a)k(Y(t,a)-X(t,a),a)
\end{equation}
(using that $J^2=-J$).
So we end up with a third version of the 
GNSB equations (\ref{GNSB1},\ref{GNSB2}),
for which:
\begin{equation}
\label{model 3}
F(x,y)=-\lambda(\hat y)k(x-\tilde y,\hat y),\;\;\;
G(x,y)=(J\mu(\hat y)k(\tilde y-x,\hat y),0).
\end{equation}
To summarize this subsection, let us just say that the GNSB equations 
(\ref{GNSB1},\ref{GNSB2})
provide a rather flexible framework to describe the interaction
between an incompressible fluid confined in a $d-$dimensional domain $D$
(each fluid parcel being labelled by its initial
position $a\in D$ having position $X(t,a)$ at time $t$),
and a set of particles (also labelled by $a$ and
of position $Y(t,a)$) moving in the
ambient space $R^d$.
The unusual feature of the resulting models is
that the interaction is pairwise ($X(t,a)$ interacts only with
$Y(t,a)$ for the same label $a$), 
except for the mediation by the pressure field $p(t,x)$ wich preserves
$a\rightarrow X(t,a)$ in the class of Lebesgue measure preserving maps of $D$
at each time $t$.

\section{The Generalized Hydrostatic Boussinesq equations}

In this final section, we investigate 
the most degenerate version of the GNSB equations 
(\ref{GNSB1},\ref{GNSB2}), where we neglect not only the inertia terms
but also the dissipative operator $K$.
Thus we are left with the strange looking system:
\begin{equation}
\label{GHB0}
F(x,y)=\nabla p,
\;\;\;\nabla\cdot v=0,
\end{equation}
\begin{equation}
\label{GHB2}
\partial_{t}y+(v\cdot\nabla)y=G(x,y),
\end{equation}
that we call Generalized Hydrostatic Boussinesq (GHB) equations
(by seing (\ref{GHB0})
as a generalization of the hydrostatic balance in Convection Theory).
Let us concentrate on the simpler case when $m=d$ and
\begin{equation}
\label{hookean}
F(x,y)=y-x
\end{equation}
(which corresponds to cables modelled by Hookean springs, according
to the mechanical interpretation of subsection \ref{mechanical}).
Thus, (\ref{GHB0}) just reads:
\begin{equation}
\label{GHB1}
y=x+\nabla p(t,x).
\end{equation}
In Lagrangian coordinates, 
the Generalized Hydrostatic Boussinesq equations
(\ref{GHB2},\ref{GHB1}) become:
\begin{equation}
\label{GHB1-lagrangian}
Y(t,a)=X(t,a)+(\nabla p)(t,X(t,a)),
\end{equation}
\begin{equation}
\label{GHB2-lagrangian}
\partial_{t}Y(t,a)=G(X(t,a),Y(t,a)),
\end{equation}
where, for all $t$, 
$a\in D\rightarrow X(t,a)\in D$ is a measure preserving map.

\subsection{Formal derivation of some optimal transport models from
the GHB equations}

We claim that several models involving optimal transport and
the Monge-Amp\`ere equation correspond to 
these GHB equations.
In particular, we consider the following generalization of the 
``semigeostrophic (SG) equations'' \cite{Ho,CNP,BB,CGP,Lo}:
\begin{equation}
\label{rho}
\partial_t \rho+\nabla\cdot (\rho w)=0,
\end{equation}
\begin{equation}
\label{w}
w(t,x)=B\nabla \varphi,
\end{equation}
\begin{equation}
\label{MA}
det(I+D^2 \varphi)=\rho,
\end{equation}
where $B$ is a $d\times d$ constant matrix and $D^2 \varphi$ is the ``Hessian''
matrix, made of all second order derivatives of $\varphi(t,x)$ with respect
to $x$.
This system, that we call Generalized Semi-Geostrophic (GSG)
equations involves the ``fully nonlinear'' 
Monge-Amp\`ere (MA) equation (\ref{MA})
which, requires, in order to be of elliptic type, the convexity condition:
\begin{equation}
\label{ellipticity}
I+D^2\varphi(t,x)>0,
\end{equation}
in the sense of symmetric matrices, for each time $t$.
The 2D SG equations \cite{Ho,CNP,BB,CGP,Lo}
just correspond to the special case when $d=2$
and $B$ is the rotation matrix of angle $\pi/2$.
The case when $B$ is just a number can be related to drift-diffusion
and Keller-Segel type models \cite{CMPS} 
for which the MA equation is
replaced by the linear Poisson equation
\begin{equation}
\label{Poisson}
\Delta \varphi=\rho-1,
\end{equation}
which can be seen as a linear approximation of the MA equation (\ref{MA}) 
as $\varphi$ is small.
The drift-diffusion case 
corresponds to 
(\ref{rho},\ref{w},\ref{Poisson}) with $B>0$.
The simplified version of the
Keller-Segel model \cite{KS} treated by J\"ager and Luckhaus \cite{JL}
corresponds to $B<0$ (with an additional diffusion term for $\rho$ in
equation (\ref{rho})).
\\
\\
Let us now show that
a solution of the
GHB equations 
(\ref{GHB1-lagrangian},\ref{GHB2-lagrangian})
corresponds to a solution of the GSG equations
(\ref{rho},\ref{w},\ref{MA},\ref{ellipticity}), in a suitable sense.
For this purpose, in order to use 
the Polar Factorization Theorem \ref{YB2}, 
we make the following $a$ $priori$ assumptions for each time $t$:
\\
A1: The map  $Y(t,\cdot)$ is non degenerate,
\\
A2: The map $x\in D\rightarrow x+\nabla p(t,x)\in R^d$
has a convex potential.

These assumptions mean that (\ref{GHB1-lagrangian}) defines the
polar factorization of $Y(t,\cdot)$ where $X(t,\cdot)$
is measure preserving and
$x\in D\rightarrow x+\nabla p(t,x)$ has a convex potential.
According to the Polar Factorization Theorem \ref{YB2}, the measure
preserving factor $X(t,\cdot)$ can be written:
$$
X(t,a)=(\nabla \Phi)(t,Y(t,a)),
$$
where $\Phi(t,x)$ is convex and Lipschitz continuous in  $x\in R^d$,
or, equivalently,
\begin{equation}
\label{XY}
X(t,a)=Y(t,a)+(\nabla \varphi)(t,Y(t,a)),
\end{equation}
where
$\varphi(t,x)=\Phi(t,x)-|x|^2/2$.
Since $Y$ is supposed to be non degenerate (by assumption A1),
there is a nonnegative Lebesgue integrable ``density field'' 
$\rho(t,x)$ such that
\begin{equation}
\label{density}
\int_{R^d}f(x)\rho(t,x)dx=\int_{D}f(Y(t,a))da,
\end{equation}
for all suitable functions $f$.
Thus:
$$
\int_{R^d}f(x+\nabla \varphi(t,x))\rho(t,x)dx
=\int_{R^d}f(Y(t,a)+(\nabla \varphi)(t,Y(t,a)))da
$$
(by definition (\ref{density}) of $\rho$)
$$
=\int_{D}f(X(t,a))da=\int_{D}f(x)dx
$$
(thanks to (\ref{XY}) and because $X(t,\cdot)$ is Lebesgue measure preserving).
So, we have obtained
\begin{equation}
\label{weak MA}
\int_{R^d}f(x+\nabla \varphi(t,x))\rho(t,x)dx=\int_{D}f(x)dx,
\end{equation}
for all compactly supported continuous function $f$.
This, combined with the assumption that $x+\nabla p(t,x)$
has a convex potential can be seen as a weak form of the Monge-Amp\`ere
equation (\ref{MA}) combined with the ellipticity condition
(\ref{ellipticity}).
Next, using (\ref{XY}), we can write (\ref{GHB2-lagrangian}) as
$$
\partial_t Y(t,a)=w(t,Y(t,a)),
$$
where $w$ is the vector field defined by:
\begin{equation}
\label{w2}
w(t,x)=G(x+\nabla \varphi(t,x),x),
\end{equation}
which is nothing but a generalization of equation (\ref{w}).
Next, we get for all smooth compactly supported function $f$ on $R^d$:
$$
\frac{d}{dt}\int_{D}f(Y(t,a))da
=\int_{D}(\nabla f)(Y(t,a))\cdot w(t,Y(t,a))da,
$$
which means, in terms of $\rho$ defined by (\ref{density}):
$$
\frac{d}{dt}\int_{R^d}f(x)\rho(t,x)dx
=\int_{R^d}\nabla f(x)\cdot w(t,x)\rho(t,x)dx,
$$
that is (\ref{rho}) in a weak sense.
So, we have fully recovered the GSG system
(\ref{rho},\ref{w},\ref{MA},\ref{ellipticity}), in a suitable weak form,
from the GHB equations (\ref{GHB1-lagrangian},\ref{GHB2-lagrangian})
under Assumption A1,A2. In addition, equation (\ref{w}) can be
replaced by the more general relation (\ref{w2}).

\subsection{A global existence theorem for the GHB equations}

We are now going to introduce a suitable concept of solutions
of (\ref{GHB1-lagrangian},\ref{GHB2-lagrangian}),
by assuming $a$ $priori$ that, in (\ref{GHB1-lagrangian}), for each
time $t$,
the map $x\in D\rightarrow x+\nabla p(t,x)$ has a convex
potential and, therefore, is the unique
rearrangement $Y^*(t,\cdot)$ of $Y(t,\cdot)$ in the
class $C$ of map with a convex potential, as in Theorem \ref{YB}.
(This corresponds to assumption A2 in the previous subsection.)
Therefore, we can just write
(\ref{GHB1-lagrangian}) as:
$$
Y(t,a)=Y^*(t,X(t,a)).
$$
From (\ref{GHB2-lagrangian}), we get (at least formally) that:
$$
\frac{d}{dt}\int_{D}f(Y(t,a))da
=\int_{D}(\nabla f)(Y(t,a))\cdot G(X(t,a),Y(t,a))da.
$$
for all compactly supported $C^1$ function $f$.
Thus,
$$
\frac{d}{dt}\int_{D}f(Y^*(t,X(t,a)))da=
$$
$$
=\int_{D}(\nabla f)(Y^*(t,X(t,a)))\cdot G(X(t,a),Y^*(t,X(t,a)))da.
$$
Now, we can factor out $X(t,a)$ and get a set of self-consistent
equations for $Y^*$, namely:
$$
\frac{d}{dt}\int_{D}f(Y^*(t,a))da
=\int_{D}(\nabla f)(Y^*(t,a))\cdot G(a,Y^*(t,a))da,
$$
without loss of information for $Y^*$.
This suggests the following concept of solution
for the GHB equations (\ref{GHB1-lagrangian},\ref{GHB2-lagrangian}):

\begin{Definition}   
\label{definition CR}
Assume that $G$ is continuous and satisfies
\begin{equation}
\label{G}
\sup_{x,y}\frac{|G(x,y)|}{1+|y|}<\infty.
\end{equation}
We say that $Y^*\in C^0([0,T],L^2(D,R^d))$ is the
``convex rearrangement'' (CR) solution 
to the GHB equations (\ref{GHB1-lagrangian},\ref{GHB2-lagrangian}), if:
\\
1) $Y^*(t,\cdot)$ belongs to the set $C$ of all maps with convex
potenial, for all $t\in [0,T]$,
\\
2) For all compactly supported $C^1$ function $f$ on $R^d$,
we have:
\begin{equation}
\label{CR solution}
\frac{d}{dt}\int_{D}f(Y^*(t,a))da
=\int_{D}(\nabla f)(Y^*(t,a))\cdot G(a,Y^*(t,a))da.
\end{equation}
\end{Definition}

This concept yields the following global existence theorem 
(without uniqueness) for all initial conditions in $L^2$:

\begin{Theorem}   
\label{existence CR}
For each initial condition
$Y^0\in L^2(D,R^d)$, there is at least one CR-solution $Y^*(t,a)$, in
the sense of Definition \ref{definition CR},
\\
such that
$Y^*(t=0,\cdot)=(Y^0)^*$. 
\\
This solution can be obtained 
as the limit in $C_t^0(L^2_a)$ as $h\rightarrow 0$
of a time discrete approximation $Y^{h}(t,a)$ defined, first
at discrete times $t=nh$, by:
\begin{equation}
\label{discrete CR}
Y^{h}(nh+h,a)=[Y^{h}(nh,a)+h\;G(a,Y^{h}(nh,a))]^*,\;\;n=0,1,2,\cdot\cdot\cdot
\end{equation}
(where $*$ denotes the rearrangement operator
as in Theorem \ref{YB})
and, then, linearly interpolated in $t$.
\end{Theorem}

\subsubsection*{Proof}

To get the existence result, it is enough to show the convergence
of the time discrete approximation $Y^h$.
First, we observe
that, $Y^h(t,\cdot)$ is valued in $C$, (the class of maps with
convex potential) for all time $t$. (This is true by definition for
discrete times $t=nh$ and preserved by linear interpolation since $C$
is a convex cone.)
Next, we deduce from (\ref{discrete CR}) and assumption (\ref{G}):
\begin{equation}
\label{bound}
\sqrt{\int_D|Y^h(nh+h,a)|^2 da}
\le hc+(1+hc)\sqrt{\int_D|Y^h(nh,a)|^2 da},
\end{equation}
for some constant $c$ depending only on $G$ and $D$.
We also get, for all compactly supported $C^1$ function $f$:
\begin{equation}
\label{consistency}
\int_D [f(Y^h(nh+h,a))-f(Y^h(nh,a))]da=
\end{equation}
$$
=h\int_D\int_0^1 (\nabla f)(Y^h(nh,a)+h\theta G(a,Y^h(nh,a)))
\cdot G(a,Y^h(nh,a))d\theta da,
$$
which can be bounded by:
\begin{equation}
\label{bound2}
hc\;\sup_x \frac{|\nabla f(x)|}{1+|x|}\int_D (1+|Y^h(nh,a)|^2)da,
\end{equation}
where $c$ depends only on $D$ and $G$.
\\
So, from (\ref{bound}), we first see that $Y^h(t,\cdot)$ is bounded in $L^2$ 
uniformly in $t\in [0,T]$ and $h\in ]0,T]$ by some constant $R$.
Therefore $Y^h$ is uniformly valued in
$C_R$ the set of maps with convex potential and $L^2$ norm bounded
by $R$, which is a compact subset of $L^2$.
\\
Next, we deduce from (\ref{consistency},\ref{bound2}) that, for each fixed
$C^1$ function $f$ such that
$$
\frac{|\nabla f(x)|}{1+|x|}<+\infty,
\;\;\;\;\;t\;\;\rightarrow\;\;\int_D f(Y^h(t,a))da
$$
is Lipschitz continuous on $[0,T]$, uniformly in $h$.
Since Theorem \ref{YB} asserts the continuity
of the map $\mu\rightarrow y^*$,
we get that $Y^h(t,\cdot)$ is uniformly
equicontinuous from $[0,T]$ to $L^2$.
Then, we deduce from the Ascoli-Arzela theorem 
that the set of all time-discrete
approximations $Y^h(t,a)$, for $0<h\le T$, is relatively compact in 
$C^0_t(L^2_a)$. Thus, there is a sequence of time steps $h$ for
which $Y^h$ converges to some limit $Y^*$ in $C^0_t(L^2_a)$, which is
necessarily valued in $C_R$, and therefore in $C$.
We also easily get (\ref{CR solution}) by letting $h$ go to zero
in (\ref{consistency}).
So, the proof of Theorem \ref{existence CR} is now complete.

\subsubsection*{Remark: continuous dependence as $d=1$}
In the very special case $d=1$,  the rearrangement 
operator $*$ is well known to
be non expansive in $L^2$:
$$
\int_D |y^*(x)-z^*(x)|^2 dx \le \int_D |y(x)-z(x)|^2 dx,
$$
for pair $(y,z)$ of $L^2$ applications from $D$ to $R$.
It follows that two CR solutions $y^*$ and $z^*$,
obtained as limits of the time-discrete approximations (\ref{discrete CR}),
with respective initial condition $y^0$ and $z^0$, must satisfy:
$$
\int_D |y^*(t,x)-z^*(t,x)|^2 dx \le \exp(ct)\int_D |y^0(x)-z^0(x)|^2 dx,
$$
where $c$ depends only on $D$ and $G$.

\section{Optimal transport and Magnetic Relaxation}

In this section, we discuss a natural ``stringy'' generalization of
the AHT equations (\ref{AHT1},\ref{AHT2}) (discussed
in the first part of the paper) and establish a link with the Arnold-Moffat
model of Magnetic Relaxation (see \cite{AK,Mo,Mo2,Sc,VMI,Ni}).

\subsection{The AHT model as a gradient flow}

Using Lagrangian coordinates, we deduce from
(\ref{AHT1},\ref{AHT2}), in the case when the
dissipation operator is $K=1$:
\begin{equation}
\label{AHT-lagrangian}
\partial_{t}X(t,a)+(\nabla p)(t,X(t,a))=y^0(a)-X(t,a)
\end{equation}
where $X(t,\cdot)$ belongs to the set $MPM(D)$ of
all (Borel) Lebesgue measure preserving maps of $D$.
(These equations can be either derived directly from 
(\ref{AHT1},\ref{AHT2}) or obtained from the GNSB equations
written in Lagrangian coordinates (\ref{GNSB-lagrangian}),
by setting $\epsilon=0$, $K=I$, $F(x,y)=y-x$ and $G(x,y)=0$.)
As explained in \cite{AHT}, in slightly different words,
the AHT equation (\ref{AHT-lagrangian}) formally corresponds to the gradient
flow of the ``energy''
\begin{equation}
\label{energy}
X\rightarrow \frac{1}{2}\int_D |X(a)-y^0(a)|^2 da,
\end{equation}
on the ``manifold'' of all $X\in MPM(D)$ for the $L^2$ metrics. 
Let us recall that, as seen in the first section, minimizing this energy
is equivalent to solve an Optimal Transport problem.
The gradient flow structure can be easily understood by considering
the standard time discretization
of such a gradient flow. Let $h>0$ be a time step
and let us denote by $X^{h}(t,a)$ the discrete approximation of $X(t,a)$
at discrete time $t=nh$, $n=0,1,2,3,\cdot\cdot\cdot.$
At $t=0$, we set $X^h(0,a)=a$ and, for $t=nh$, $n=1,2,3,\cdot\cdot\cdot$,
we require $X^h(t,\cdot)$ to be
a minimizer among all $X\in MPM(D)$
of the following functional:
\begin{equation}
\label{discrete AHT}
X\rightarrow \int_D \frac{|X(a)-X^{h}(nh-h,a)|^2}{2h}
+\frac{1}{2}\int_D |X(a)-y^0(a)|^2 da,
\end{equation}
or, equivalently,
\begin{equation}
\label{discrete AHT2}
\int_D |X(a)-\frac{X^{h}(nh-h,a)+y^0(a)h}{1+h}|^2 da
\end{equation}
(after rearranging the squares).
\\
Thus, assuming $a$ $priori$
that $X^{h}(nh-h,\cdot)+y^0 h$ is non degenerate and using
Theorem \ref{YB2}, this exactly means:
$$
X^{h}(nh-h,\cdot)+y^0 h=y^*\circ X^{h}(nh,\cdot)
$$
where $y^*$ is a map with a convex potential.
If we write this map as
$$
x\rightarrow x+h(\nabla p^h)(nh,x),
$$
we get
$$
X^{h}(nh-h,\cdot)+y^0 h=X^{h}(nh,\cdot)+h(\nabla p^h)(nh,X^{h}(nh,\cdot)),
$$
which can be seen just as a finite difference approximation of 
equation (\ref{AHT-lagrangian}) as $h\rightarrow 0$.
This is enough to interpret, at least formally, equation
(\ref{AHT-lagrangian}) as the $L^2$ gradient flow of energy (\ref{energy})
on the ``manifold'' $MPM(D)$.

\subsection{A stringy generalization of the AHT model}

This analysis suggests
the possibility of more complex models based on similar ideas.
A rather natural idea amounts to consider, instead of the ``manifold'' 
$MPM(D)$, the ``manifold'' of ``strings'' valued in $MPM(D)$:
$s\in [0,1]\rightarrow X(\cdot,s)\in MPM(D)$,
with fixed end values, say 
\begin{equation}
\label{end points}
X(t,a,s=0)=X^-(t,a),\;\;\;X(t,a,s=1)=X^+(t,a).
\end{equation}
Then, we may think of the $L^2$ gradient flow of the following string
energy:
\begin{equation}
\label{energy2}
\frac{1}{2}\int_0^1\int_D |\partial_s X(a,s)|^2 da\;ds.
\end{equation}
We claim that the resulting equation read:
\begin{equation}
\label{MR-lagrangian}
\partial_t X(t,a,s)=\partial^2_{ss}X(t,a,s)
+(\nabla p)(t,X(t,a,s),s),
\end{equation}
where $X(t,\cdot,s)$ is valued in $MPM(D)$ and end point conditions
(\ref{end points}) are enforced.
To get this system, as we did before for the AHT model,
we define a time discrete approximation 
$X^h(t,a,s)$,
by setting $X^h(0,a,s)=a$ and asking, for $t=nh$, $n=1,2,3,\cdot\cdot\cdot$,
$X^h(t,\cdot,\cdot)$ to be
a minimizer among all curves $s\in [0,1]\rightarrow X(\cdot,s)\in MPM(D)$
of the functional:
\begin{equation}
\label{discrete MR}
\int_D\int_0^1  \frac{|X(a,s)-X^{h}(nh-h,a,s)|^2}{2h}da\;ds
+\int_D\int_0^1  \frac{|\partial_s X(a,s)|^2}{2}da\;ds.
\end{equation}
The formal optimality condition reads:
$$
\frac{X^{h}(nh,a,s)-X^{h}(nh-h,a,s)}{h}=
$$
$$
=\partial^2_{ss}X^{h}(nh,a,s)
+(\nabla p^h)(nh,X^{h}(nh,a,s),s),
$$
for some scalar function $p^h$.
So, we formally obtain, as $h\rightarrow 0$, the desired
equation (\ref{MR-lagrangian}).
Equation (\ref{MR-lagrangian}) has an interesting interpretation,
obtained by assuming $a$ $priori$ that $a\in D\rightarrow X(t,a,s)$
is a smooth orientation and measuring preserving diffeomorphism of $D$
for each $(t,s)$.
Then we introduce, for each $(t,s)$, 
two divergence free vector fields parallel to the boundary $\partial D$,
namely $v(t,x,s)\in R^d$ and $b(t,x,s)\in R^d$, defined by:
\begin{equation}
\label{def v b}
v(t,X(t,a,s),s)=\partial_t X(t,a,s),
\;\;\;b(t,X(t,a,s),s)=\partial_s X(t,a,s).
\end{equation}
Then, we get from (\ref{MR-lagrangian}):
\begin{equation}
\label{MR1}
v=\partial_s b+(b\cdot \nabla)b+\nabla p,
\end{equation}
while, from (\ref{def v b}), we get the compatibility condition
\begin{equation}
\label{MR2}
\partial_t b+(v\cdot \nabla)b=(b\cdot \nabla)v+\partial_s v
\end{equation}
(by writting $\partial^2_{ts}X=\partial^2_{st}X$),
to be added to the divergence free constraints
\begin{equation}
\label{MR3}
\nabla\cdot v=\nabla\cdot b=0,
\;\;\;v//\partial D,\;\;\;b//\partial D,
\end{equation}
and the boundary conditions at $s=0$ and $s=1$ induced by
(\ref{end points}), namely:
\begin{equation}
\label{MR4}
v(t,x,s=0)=v^-(t,x),\;\;\;v(t,x,s=1)=v^+(t,x),
\end{equation}
where $v^+$ and $v^-$ are prescribed.
When the fields $v$ and $b$ do not depend on $s$, we get the 
Magnetic Relaxation model discussed by Moffatt in \cite{Mo}
(see also \cite{AK,Mo2,Sc,VMI,Ni}).
As $t\rightarrow+\infty$, we expect, at least for a large class
of initial conditions, the solution 
of equations (\ref{MR1},\ref{MR2},\ref{MR3})
to converge toward
an equilibrium, for which $v=0$ and $b=b(x,s)$, $p=p(x,s)$ are
solutions to the Euler equations \cite{AK,MP}
($s$ acting as the time variable
and $b$ as the velocity field):
\begin{equation}
\label{Euler1}
\partial_s b+(b\cdot \nabla)b+\nabla p=0,
\end{equation}
\begin{equation}
\label{Euler2}
\nabla\cdot b=0,
\;\;\;b//\partial D.
\end{equation}
Of course, we are far from being able
to provide any rigorous proof of this conjecture.

\subsection{The ``cross-Burgers'' equation}

In the case when $D$ is the unit ball,
The Magnetic Relaxation equations 
(\ref{MR1},\ref{MR2},\ref{MR3})
admit special solutions $(b,v,\nabla p)$
which are linear in $x$:
\begin{equation}
\label{ansatz}
b(t,x,s)=B(t,s)x,\;\;\;
v(t,x,s)=V(t,s)x,\;\;\;
\nabla p(t,x,s)=G(t,s)x,
\end{equation}
where $B$, $V$ are skew-symmetric matrices, while $G$ is a symmetric
matrix, all depending only on $(t,s)$.
(Notice that the fields $b$ and $v$ are automatically parallel
to the boundary $\partial D$ since $D$ is the unit ball.)
The resulting equations for $B$, $V$ and $G$ are:
\begin{equation}
\label{MR1-toy}
V=\partial_s B+B^2+G,
\end{equation}
\begin{equation}
\label{MR2-toy}
\partial_t B+[V,B]=\partial_s V.
\end{equation}
Since $B^2$ is a symmetric matrix, equation (\ref{MR1-toy})
reduces to:
$$
V=\partial_s B.
$$
Thus, we get a single equation for B:
\begin{equation}
\label{bracket burgers}
\partial_t B+[\partial_s B,B]=\partial^2_{ss}B,
\end{equation}
where $[A,B]$ denotes the skew product $AB-BA$.
In the
special case $d=3$, $B$ can be identified as a 3-vector and
$[\cdot,\cdot]$ as the cross product $\times$ in $R^3$, which leads
to:
\begin{equation}
\label{cross burgers}
\partial_t B+\partial_s B\times B=\partial^2_{ss}B.
\end{equation}
that we could call the ``cross-Burgers'' equation.
This equation admits interesting special solutions, such as:
$$
B(t,s)=(\alpha(t)\cos\;s\;,\alpha(t)\sin\;s\;,\beta(t)-1)
$$
where $\alpha\ge 0$ and $\beta$ are solutions to:
$$
\frac{d\alpha}{dt}=-\beta\alpha,
\;\;\;\frac{d\beta}{dt}=\alpha^2,
$$
or, equivalently, 
$$
\frac{d^2\lambda}{dt^2}+\exp(2\lambda)=0,
$$
where $\lambda=\log(\alpha)$.

\section*{Acknowledgment}
This work originated 
while the author was visiting 
Universit\"at Bonn (Institut f\"ur Angewandte Mathematik)
and continued at Universit\"at Wien (Wolfgang Pauli Institut).
This work is also part of the research 
made in the LRC CEA-Cadarache/CNRS-Universit\'e de Nice
and the ANR OTARIE project (ANR BLAN07-2-183172).
We thank Adrien Blanchet, Marco Di Francesco, Francis Filbet
and Fran\c cois Gallaire for fruitful discussions.


\begin{thebibliography}{A}
\bibitem[Am]{Am} L. Ambrosio, 
{\it Transport equation and Cauchy problem for $BV$ vector fields.}
{\sl Invent. Math. 158 (2004) 227-260.}
\bibitem[AK]{AK} V. Arnold, B. Khesin, 
{\it Topological methods in hydrodynamics,}
{\sl Applied Mathematical Sciences, 125. Springer, 1998.}
\bibitem[AHT]{AHT} S. Angenent, S. Haker, A. Tannenbaum, 
{\it Minimizing flows for the Monge-Kantorovich problem,}
{\sl  SIAM J. Math. Anal.  35  (2003) 61-97.}
\bibitem[BB]{BB} J.-D. Benamou, Y. Brenier,
{\it Weak existence for the semigeostrophic equations formulated as a coupled 
Monge-Amp\`ere/transport problem,}
{\sl SIAM J. Appl. Math. 58 (1998) 1450-1461.}
\bibitem[Br]{Br} Y. Brenier,
{\it Polar factorization and monotone rearrangement 
of vector-valued functions,} 
{\sl Comm. Pure Appl. Math., 64 (1991) 375-417.}
\bibitem[Ca]{Ca} L. Caffarelli, 
{\it Boundary regularity of maps with convex potentials,}
{\sl Comm. Pure Appl. Math. 45 (1992) 1141-1151.}
\bibitem[Ch]{Ch} D. Chae,
{\it Global regularity for the 2D Boussinesq equations 
with partial viscosity terms,}
{\sl  Adv. Math.  203  (2006) 497-513.}
\bibitem[CMPS]{CMPS} F. Chalub, P. Markowich, B. Perthame, C. Schmeiser,
{\it Kinetic models for chemotaxis and their drift-diffusion limits,}
{\sl Monatsh. Math. 142 (2004) 123-141.}
\bibitem[CGP]{CGP} M. Cullen, W. Gangbo, G. Pisante,
{\it The semigeostrophic equations discretized in reference and dual 
variables,}
{\sl  Arch. Ration. Mech. Anal.  185  (2007) 341-363.}
\bibitem[CNP]{CNP} M. Cullen, J. Norbury, J. Purser, 
{\it Generalised Lagrangian
solutions for atmospheric and oceanic flows,}
{\sl SIAM J. Appl. Math. 51 (1991) 20-31.}
\bibitem[DL]{DL} R. Di Perna,, P.-L. Lions,
{\it Ordinary differential equations, transport theory and Sobolev spaces,}
{\sl Invent. Math. 98 (1989) 511-547.}
\bibitem[DOR]{DOR} C. Doering, F. Otto, M. Reznikoff, 
{\it Bounds on vertical heat transport for infinite-Prandtl-number 
Rayleigh-B\'enard convection,}
{\sl  J. Fluid Mech.  560  (2006) 229-241.}
\bibitem[FMMS]{FMMS} U. Frisch, S. Matarrese, R. Mohayaee, A. Sobolevski. 
{\it A reconstruction of the initial conditions of the Universe by optimal 
mass reconstruction,}
{\sl Nature, 417 (2002) 260-262.}
\bibitem[HZTA]{HZTA} S. Haker, L. Zhu, A. Tannenbaum, and S. Angenent,
{\it Optimal Mass Transport for Registration and Warping,}
{\sl International Journal on Computer Vision, 60(3) (2004) 225-240.}
\bibitem[Ho]{Ho} B. Hoskins, 
{\it The mathematical theory of frontogenesis,}
{\sl Annual review
of fluid mechanics, Vol. 14, pp. 131-151, Palo Alto, 1982.}
\bibitem[HL]{HL} T. Hou, C. Li,
{\it Global well-posedness of the viscous Boussinesq equations,}
{\sl Discrete Contin. Dyn. Syst. 12 (2005) 1-12.}
\bibitem[ID]{ID} FP Incropere, DP DeWitt,
{\it Heat and Mass Transfer,}
{\sl John Wiley and Sons, New York, 1996}
\bibitem[JL]{JL} W. J\"ager, S. Luckhaus, 
{\it On explosions of solutions to a system of partial differential equations 
modelling chemotaxis,}
{\sl Trans. Amer. Math. Soc. 329 (1992) 819-824.}
\bibitem[KS]{KS} E. Keller, L. Segel
{\it Model for chemotaxis,}
{\sl J Theor Biol 30 (1971) 225-234.}
\bibitem[Li]{Li} P.-L. Lions,
{\it Mathematical topics in fluid mechanics. Vol. 1. Incompressible models,}
{\sl Oxford Lecture Series in Mathematics
and its Applications, Oxford University Press, New York, 1996.}
\bibitem[Lo]{Lo} G. Loeper,
{\it A fully nonlinear version of the incompressible Euler equations: 
the semigeostrophic system,}
{\sl SIAM J. Math. Anal.  38  (2006) 795-823.}
\bibitem[Mc]{Mc} R. McCann, 
{\it Polar factorization of maps on Riemannian manifolds,}
{\sl Geom. Funct. Anal. 11 (2001) 589-608.}
\bibitem[Ma]{Ma} A. Majda,
{\it Introduction to PDEs and Waves for the Atmosphere and Ocean,}
{\sl AMS and CIMS, 2000.}
\bibitem[MP]{MP} C.Marchioro, M.Pulvirenti,
{\it Mathematical theory of incompressible nonviscous fluids,}
{\sl Springer, New York, 1994.}
\bibitem[Mo]{Mo} H.K. Moffatt,
{\it Magnetostatic equilibria and analogous Euler flows of 
arbitrarily complex topology. I. Fundamentals,}
{\sl J. Fluid Mech. 159 (1985) 359-378.}
\bibitem[Mo2]{Mo2} H. Moffatt,
{\it Relaxation under topological constraints,}
{\sl Topological aspects of the dynamics of fluids and plasmas,
NATO Adv. Sci. Inst. Ser. E Appl. Sci., 218, Kluwer, 1992.}
\bibitem[NPS]{NPS} J. Nieto, F. Poupaud, J. Soler, 
{\it High-field limit for the Vlasov-Poisson-Fokker-Planck system,}
{\sl  Arch. Ration. Mech. Anal.  158  (2001) 29-59.}
\bibitem[Ni]{Ni} T. Nishiyama,
{\it Magnetohydrodynamic approaches to measure-valued solutions of the 
two-dimensional stationary Euler equations,}
{\sl  Bull. Inst. Math. Acad. Sin. (N.S.)  2  (2007) 139-154.}
\bibitem[Pe]{Pe} J.Pedlosky,
{\it Geophysical fluid dynamics,}
{\sl Springer, New York, 1979.}
\bibitem[Sc]{Sc} M. Schonbek, 
{\it Decay of Solutions to non-oscillating Magneto Hydrodynamics equations,}
{\sl Theory of the Navier-Stokes equations,  179-184, Ser. Adv. Math. Appl. 
Sci., 47, World Sci., 1998.}
\bibitem[Vi]{Vi} C. Villani,
{\it Topics in optimal transportation,}
{\sl Graduate Studies in Mathematics, 58, AMS, Providence, 2003.}
\bibitem[VMI]{VMI} V.A. Vladimirov, H.K. Moffatt, K.I. Ilin,
{\it On general transformations and variational principles for the 
magnetohydrodynamics of ideal fluids. IV,}
{\sl J. Fluid Mech. 390 (1999) 127-150.}
\bibitem[Wa]{Wa} X. Wang,
{\it Infinite Prandtl number limit of Rayleigh-B\'enard convection,}
{\sl Comm. Pure Appl. Math. 57 (2004) 1265-1282.}
\end{thebibliography}
\end{document}